\RequirePackage{fix-cm}
\documentclass{svjour3}                     % onecolumn (standard format)
\smartqed  % flush right qed marks, e.g. at end of proof
\def\makeheadbox{\relax} % for submission comment it
\usepackage{basic}
\usepackage{amsmath}
\usepackage{amsfonts}
\usepackage{amsmath}
\usepackage{times}
\usepackage{dsfont}
\begin{document}

\title{Optimal rates for F-score binary classification}
% \thanks{This work is supported by the Labex B\'ezout of the Universit\'e Paris-Est.}
% \subtitle{$\text{F}_1$-score binary classification}

\titlerunning{F-score binary classification}        % if too long for running head

\author{Evgenii Chzhen}


\institute{LAMA, Universit\'e Paris-Est \at
              Cit\'e Descartes\\
              5 boulevard Descartes\\
              77454 Marne-la-Vall\'ee cedex 2\\
              \email{evgenii.chzhen@univ-paris-est.fr}
}
% \date{Received: date / Accepted: date}  % uncomment for submission
\date{}
\maketitle

\begin{abstract}
  We study the minimax settings of binary classification with $\text{F}$-score under the $\beta$-smoothness assumptions on the regression function $\eta(x) = \Prob(Y = 1| X = x)$ for $x \in \bbR^d$.
    We propose a classification procedure which under the $\alpha$-margin assumption achieves the rate $\bigO(n^{-(1 + \alpha)\beta / (2\beta + d)})$ for the excess $\text{F}$-score.
    In this context, the Bayes optimal classifier for the $\text{F}$-score can be obtained by thresholding the aforementioned regression function $\eta$ on some level $\theta^*$ to be estimated.
    The proposed procedure is performed in a semi-supervised manner, that is, for the estimation of the regression function we use a labeled dataset of size $n \in \bbN$ and for the estimation of the optimal threshold $\theta^*$ we use an unlabeled dataset of size $N \in \bbN$.
    Interestingly, the value of $N \in \bbN$ does not affect the rate of convergence, which indicates that it is ``harder'' to estimate the regression function $\eta$ than the optimal threshold $\theta^*$.
    This further implies that the binary classification with $\text{F}$-score behaves similarly to the standard settings of binary classification.
    Finally, we show that the rates achieved by the proposed procedure are optimal in the minimax sense up to a constant factor.
\end{abstract}

%%%%%%%%%%%%%%%%%%%%%%%%%%%%%%%%%%%%%%%%%%%%%%%%%%%%%%%%%%%%%%%%%%%%%%%%%%%%%%%
\section{Introduction}
%%%%%%%%%%%%%%%%%%%%%%%%%%%%%%%%%%%%%%%%%%%%%%%%%%%%%%%%%%%%%%%%%%%%%%%%%%%%%%%
The problem of binary classification is among the most basic and well-studied problems in statistics and machine learning~\cite{Vapnik98,Yang99,Bartlett_Mendelson02,Audibert04,Massart_Nedelec06,Audibert_Tsybakov07}.
Until very recently, theoretical guarantees were almost exclusively formulated in terms of the probability of miss-classification (a.k.a accuracy) as the measure of the risk.
This choice of the risk is practically suitable in the case of the ``well-balanced'' distributions and datasets, that is, the probabilities to observe both classes are similar.

Once this assumption fails to be satisfied, classifiers based on the accuracy might perform poorly in practice.
One possible approach to treat such a situation is to modify the measure to be optimized in an appropriate way.
A popular choice of such measure is the F-score, whose roots can be tracked back to the information retrieval literature~\cite{Rijsbergen74,Lewis95}.
From the statistical point of view there are two alternative approaches~\cite{Ye_Chai_Lee_Chieu12,Dembczynski_Kotlowski_Koyejo_Natarajan17} to the theoretical treatment of the F-score: Population Utility (PU) and Expected Test Utility (ETU).
In this work we follow the PU approach which, as noted in~\cite{Dembczynski_Kotlowski_Koyejo_Natarajan17}, has stronger roots in classical statistics.
Our goal is to provide minimax analysis of the binary classification with F-score under non-parametric assumptions.

%%%%%%%%%%%%%%%%%%%%%%%%%%%%%%%%%%%%%%%%%%%%%%%%%%%%%%%%%%%%%%%%%%%%%%%%%%%%%%%
\section{The problem formulation}
%%%%%%%%%%%%%%%%%%%%%%%%%%%%%%%%%%%%%%%%%%%%%%%%%%%%%%%%%%%%%%%%%%%%%%%%%%%%%%%
We first introduce some notation that is used throughout this work.
For any two real numbers $a, b \in \bbR$ we denote by $a \wedge b$ (resp. $a \vee b$) the minimum (resp the maximum) between $a$ and $b$.
The standard Euclidean norm in $\bbR^d$ is denoted by $\norm{\cdot}_2$ and a ball centered at $x \in \bbR^d$ of radius $r$ is denoted by $\class{B}(x, r)$.
For positive real valued sequences $a_n, b_n : \bbN \mapsto \bbR_+$ we say that $a_n = \bigO(b_n)$ if there exists some positive constant $M > 0$ such that for all $n \in \bbN$ it holds that $a_n / b_n \leq M$.
We consider a random couple $(X, Y)$ taking values in $\bbR^d \times \{0, 1\}$ with joint distribution $\Prob$.
The vector $X \in \bbR^d$ is the feature vector and the binary variable $Y \in \{0, 1\}$ is the label, in what follows we assume that $\Prob(Y = 1) \neq 0$.
We denote by $\Prob_X$ the marginal distribution of the feature vector $X \in \bbR^d$ and by $\eta(X) \coloneqq \Prob(Y = 1 | X)$ the regression function.
A classifier is any measurable function $g : \bbR^d \mapsto \{0, 1\}$ and the set of all such functions is denoted by $\class{G}$.

We assume that we have access to two datasets: the first dataset $\data_n = \{(X_i, Y_i)\}_{i = 1}^n$ consists of $n \in \bbN$ \iid copies of $(X, Y) \sim \Prob$; and the second dataset $\data_N = \{X_i\}_{i = n+1}^{n + N}$ consists of $N \in \bbN$ independent copies of $X \sim \Prob_X$.
Denote by $\Prob^{\otimes n}$ and $\Prob_X^{\otimes N}$ the distributions of $\data_n$ and $\data_N$ respectively.
Moreover, we denote by $\Exp_{(\data_n, \data_N)}$ the expectation with respect to the distribution of $(\data_n, \data_N)$, that is, with respect to $\Prob^{\otimes n} \otimes \Prob_X^{\otimes N}$ on the space $\left(\bbR^d \times \{0, 1\}\right)^{n} \times \left(\bbR^d\right)^{N}$.
We additionally assume that the size of the unlabeled dataset is not smaller that the size of the labeled dataset\footnote{Note that one can always satisfy this assumption by augmenting $\data_N$ using a portion of $\data_n$ and erasing labels. Typically, in practice it is easier to gather the unlabeled data then labeled, that is why this assumption is rather a formality.}, that is, $N \geq n$.
For a given classifier $g: \bbR^d \mapsto \{0, 1\}$ we define its $\text{F}_b$-score\footnote{We decided to divide the classical definition of the $\text{F}_b$-score by the factor $1 + b^2$ to simplify the notation, thus, it is sufficient to multiply the obtained results by $1 + b^2$, to recover the results on the classical definition of the $\text{F}_b$-score.} for any $b>0$ by
\begin{align*}
    F_b(g) \coloneqq \frac{\Prob(Y = 1, g(X) = 1)}{b^2\Prob(Y = 1) + \Prob(g(X) = 1)}\enspace.
\end{align*}
A Bayes-optimal classifier $g^*: \bbR^d \mapsto \{0, 1\}$ is any classifier that maximizes the $\text{F}$-score over all classifiers $\class{G}$, that is,
\begin{align*}
    g^* \in \argmax_{g \in \class{G}}F_b(g)\enspace.
\end{align*}
It was established by \cite{Zhao_Edakunni_Pocock_Brown13} that a maximizer of the $\text{F}_1$-score can be obtained by comparing the regression function $\eta(X)$ with a threshold $\theta^* \in [0, 1]$.
Importantly, this threshold depends explicitly on the distribution $\Prob$ and can be obtained as unique root of
\begin{align*}
    \theta \mapsto \theta \Prob(Y = 1) - \Exp(\eta(X) - \theta)_+\enspace.
\end{align*}
One of the contributions of this work is extension of the result of~\cite[Section 6]{Zhao_Edakunni_Pocock_Brown13} for an arbitrary value of $b > 0$.
\begin{theorem}
    \label{thm:bayes_optimal}
    A Bayes-optimal classifier $g^*$ can be obtained point-wise for all $x \in \bbR^d$ as
    \begin{align}
        g^*(x) = \ind{\eta(x) > \theta^*}\enspace,
    \end{align}
    where $\theta^* \in [0, 1]$ is a threshold which satisfies
    \begin{align*}
        b^2\theta^* \Prob(Y = 1) = \Exp(\eta(X) - \theta^*)_+\enspace.
    \end{align*}
    Moreover, the classifier $g^*$ satisfies $F_b(g^*) = \theta^*$.
\end{theorem}
The proof can be found in Appendix~\ref{app:bayes_excess}
Notice that if the optimal threshold $\theta^* \in [0, 1]$ is known a priori, the problem of binary classification with the $\text{F}$-score is no harder than the standard settings of binary classification with the accuracy as the measure of performance.
As the threshold $\theta^* \in [0, 1]$ depends on the distribution $\Prob$, it could be estimated using data.
Theorem~\ref{thm:bayes_optimal} allows to obtain a trivial upper bound on the threshold $\theta^*$, indeed, since $\theta^* = F_b(g^*)$ and for any classifier $g \in \class{G}$ the $\text{F}_b$-score is upper bounded by $1/(1 + b^2)$ we have $\theta^* \in [0, 1/(1 + b^2)]$.

For any classifier $g: \bbR^d \mapsto \{0, 1\}$ we define its excess score as
\begin{align*}
    \excess_b(g) \coloneqq F_b(g^*) - F_b(g),\quad\text{(excess score)}\enspace.
\end{align*}
The excess score is the central object of our analysis and one of our goals is to provide an estimator whose excess score is as small as possible.
Using Theorem~\ref{thm:bayes_optimal} we can show that the excess score of any classifier $g: \bbR^d \mapsto \{0, 1\}$ can be written in a simple form.
\begin{lemma}
    \label{lem:excess_score}
    Let $g: \bbR^d \mapsto \{0, 1\}$ be any classifier and assume that $\Prob(Y = 1) \neq 0$, then
    \begin{align*}
        \excess_b(g) = \frac{\Exp\abs{\eta(X) - \theta^*}\ind{g^*(X) \neq g(X)}}{b^2\Prob(Y = 1) + \Prob(g(X) = 1)}\enspace.
    \end{align*}
\end{lemma}
In general the Bayes optimal rule is not unique, Theorem~\ref{thm:bayes_optimal} only states that one of the optimal classifiers has the form described by its statement.
Even though, the function $\theta \mapsto b^2\theta \Prob(Y = 1) - \Exp(\eta(X) - \theta)_+$ has unique root (see Appendix~\ref{app:upper} for the proof), other thresholds may result in the same Bayes rule.
Indeed, consider a simple example with $\eta(x) \equiv 1/2$, $b = 1$, then it is easy to see that the solution $\theta^*$ of $\theta / 2 = (1/2 - \theta)_+$ is exactly $1/3$,
and every Bayes optimal classifier predicts one almost surely.
Clearly, any threshold $\theta \in [0, 1/2)$ of the regression function $\eta$ results in the same classifier.
Importantly, Lemma~\ref{lem:excess_score} and the equality $\argmax_{g \in \class{G}} F_1(g) = \theta^*$ are valid \emph{only} for the threshold $\theta^* = 1/3$.
In this work, we shall always refer to $\theta^*$ being the solution of $b^2\theta \Prob(Y = 1) = \Exp(\eta(X) - \theta)_+$ and we call this threshold as the \emph{optimal} threshold.
\begin{remark}
    For the rest of the paper, we focus our attention only on the value $b = 1$ to simplify the presentation.
    It will be clear from our arguments that the generalization of the theoretical results of the paper to an arbitrary value $b > 0$ follows straightforwardly from our analysis.
\end{remark}
Interestingly, the results above demonstrate that the problem of binary classification with F-score has a lot in common with the standard settings.
Indeed, in both cases the Bayes optimal classifier is obtained via thresholding of the regression function and the expression for the excess risk is also similar.
Consequently, in this work we address the following questions
\begin{itemize}
    \item[Q1.:] Is the problem of binary classification with F-score harder than its more known counterpart?
    In particular, can the minimax analysis of~\cite{Audibert_Tsybakov07} be extended to these settings and what is an optimal algorithm?
    \item[Q2.:] In view of recent results of~\cite{Chzhen_Denis_Hebiri19}, we wonder if the introduction of unlabeled dataset can improve classification algorithms in the context of F-score.
\end{itemize}

Let us point out, that Lemma~\ref{lem:excess_score} is crucial for our analysis as it allows to adapt the scheme provided by~\cite{Audibert_Tsybakov07} for the standard setting of the binary classification.
Nevertheless, as the threshold $\theta^* \in [0, 1]$ is unknown beforehand, this machinery cannot be applied in a straightforward way and some effort is required.
In this work, we pose similar assumptions on the distribution $\Prob$ to the ones used in~\cite{Audibert_Tsybakov07}.
\begin{assumption}[$\alpha$-margin assumption]
    \label{ass:margin_assumption}
    We say that the distribution $\Prob$ of the pair $(X, Y) \in \bbR^d \times \{0, 1\}$ satisfies the $\alpha$-margin assumption if there exist constants $C_0 > 0$, $\delta_0 \in (0, 1/12]$ and $\alpha > 0$ such that for every positive $\delta \leq \delta_0$ we have
    \begin{align*}
        \Prob_X(0 < \abs{\eta(X) - \theta^*} \leq \delta) \leq C_0 \delta^\alpha\enspace.
    \end{align*}
\end{assumption}
    The case of ``$\alpha = \infty$'' is understood in the following manner~\cite{Massart_Nedelec06}:
    there exists a constant $\delta_0 \in (0, 1]$ such that
    \begin{align*}
        \Prob_X(0 < \abs{\eta(X) - \theta^*} \leq \delta_0) = 0\enspace,
    \end{align*}
    typically this is the most advantageous situation for the binary classification, as the regression function $\eta$ is separated from the optimal threshold $\theta^*$.
    Assumption~\ref{ass:margin_assumption} specifies the concentration rate of the regression function $\eta$ around the optimal threshold $\theta^*$.
    Notice, if Assumption~\ref{ass:margin_assumption} is satisfied, it holds that for all $\delta > 0$
    \begin{align*}
        \Prob_X(0 < \abs{\eta(X) - \theta^*} \leq \delta) \leq c_0 \delta^\alpha\enspace,
    \end{align*}
    where $c_0 = C_0 \vee \delta_0^{-\alpha}$.
This condition is tightly related to the rate of convergence in the case of the binary classification~\cite{Audibert_Tsybakov07,Massart_Nedelec06}.
The classification algorithm that is proposed in this work is based on a direct estimation of the regression function $\eta$ and the optimal threshold $\theta^*$.
% To this end, we assume that the regression function $\eta: \bbR^d \mapsto [0, 1]$ can be estimated accurately.

In the sequel, we consider the case of non-parametric estimation, that is we assume that the regression function $\eta: \bbR^d \mapsto \{0, 1\}$ lies in some class of $\beta$-smooth functions and the marginal density $\Prob_X$ of $X \in \bbR^d$ admits density \wrt to the Lebesgue measure supported on a well-behaved compact set and uniformly lower- and upper bounded.
The exact formal description of these assumptions is given in Section~\ref{sec:lower_bound}, where we prove optimality of our rates.
As for now, it is sufficient to assume that there exists a good estimator $\heta$ based on the labeled set $\data_n$ of the regression function $\eta$.
\begin{assumption}[Existence of estimator]
    \label{ass:exponential_concentration_estimator}
    There exists an estimator $\heta$ based on $\data_n$ which satisfies for all $t > 0$
    \begin{align*}
        \Prob^{\otimes n}(\abs{\heta(x) - \eta(x)} \geq t) \leq C_1\exp(-C_2a_n t^2) \text{ a.s. } \Prob_X\enspace,
    \end{align*}
    for some universal constants $C_1, C_2 > 0$ and an increasing sequence $a_n: \bbN \mapsto \bbR_+$.
\end{assumption}
For instance, in the case of $\beta$-smooth regression function\footnote{Typically, one also need to assume that the marginal distribution $\Prob_X$ is well behaved, see Section~\ref{sec:lower_bound}.} $\eta: \bbR^d \mapsto [0, 1]$, a typical non-parametric rate is $a_n = n^{{2\beta}/{(2\beta + d)}}$ and it can be achieved by the local polynomial estimator, see \cite[Theorem 3.2]{Audibert_Tsybakov07}.
Finally, in this work we assume that the probability $\Prob(Y = 1)$ is lower bounded by some constant which can be arbitrary small but fixed.
\begin{assumption}[Lower bounded $\Prob(Y = 1)$]
    \label{ass:lower_proba_1}
    We assume that there exists a positive constant $p$ such that $ p \leq \Prob(Y = 1)$.
\end{assumption}
It is assumed that the constants $C_0, C_1, C_2, p$ are independent of both $n, N \in \bbN$, however these constants can depend on the dimension of the problem $d$, on the value of $\alpha > 0$ as well as on each other.
The values of the constants $C_0, C_1, C_2, p$ are not going to impact the rates of convergence, though they might and will enter as numerical constants in front of the rate.
In contrast, the value of $\alpha$ in the margin assumption will explicitly appear in the obtained rates.
%%%%%%%%%%%%%%%%%%%%%%%%%%%%%%%%%%%%%%%%%%%%%%%%%%%%%%%%%%%%%%%%%%%%%%%%%%%%%%%
\section{Related works and contributions}
%%%%%%%%%%%%%%%%%%%%%%%%%%%%%%%%%%%%%%%%%%%%%%%%%%%%%%%%%%%%%%%%%%%%%%%%%%%%%%%
Literature on the binary classification with $\text{F}$-score is rather broad, it spans both applied and theoretical studies of the problem.
It should be noted that our work falls into the Population Utility (PU) approach~\cite{Dembczynski_Kotlowski_Koyejo_Natarajan17}, that is, the expectation is taken in the numerator and the denominator of the $\text{F}$-score simultaneously.
This approach should not be confused with the Expected Test Utility (ETU) approach, for which a non-asymptotic behavior can differ significantly.
We refer the reader to~\cite{Dembczynski_Kotlowski_Koyejo_Natarajan17,Ye_Chai_Lee_Chieu12} where the PU vs. ETU tale is discussed in depth and their asymptotic equivalency is established.
Let us mention that, the asymptotic statistical theory of the binary classification with $\text{F}$-score has been studied in the prior literature~\cite{Koyejo_Natarajan_Ravikumar_Dhillon14,Narasimhan_Vaish_Agarwal14,Menon_Narasimhan_Agarwal_Chawla13,Ye_Chai_Lee_Chieu12}.
Bellow, we summarize our contributions and highlight the improvements with respect to the previous results on the non-asymptotic analysis of the binary classification with $\text{F}$-score.
\begin{itemize}
    \item We propose a two-step estimator, which first estimates the regression function $\eta$ and then the optimal threshold $\theta^*$.
    This type of two-step estimators, which involve an explicit thresholds tuning, are well-known in the literature and demonstrate a promising empirical performance~\cite{Koyejo_Natarajan_Ravikumar_Dhillon14,Keerthi_Sindhwani_Chapelle07}.
    An important novelty introduced here is the semi-supervised nature of the procedure which can exploit the unlabeled data.
    It is already a well established fact that the semi-supervised methods might~\cite{Singh_Nowak_Zhu09} or not~\cite{Rigollet07} improve supervised estimation from statistical point of view.
    However, let us point out, that from the practical point of view, typically the most expensive part of the data gathering process is the correct labeling.
    Thus, one may assume that the unlabeled dataset $\data_N$ is always available in reality and the settings $N \gg n$ are satisfied.
    Our analysis implies that in the setting of binary classification with F-score the semi-supervised techniques are not superior to the supervised ones.
    % Similar conclusion were established in~\cite{Rigollet07} in the context of classification under the cluster assumption.
    In contrast, in~\cite{Chzhen_Denis_Hebiri19} the authors showed that in the context of confidence set classification semi-supervised classifiers might outperform it supervised counterparts.
    \item From the theoretical point of view, the most relevant reference is a recent work of~\cite{Yan_Koyejo_Zhong_Ravikumar18}, where the authors studied a rather broad class of performance measures for the problem of binary classification, namely Karmic measures.
    This class includes the $\text{F}$-score, considered in the present manuscript.
    Under similar, though stronger assumptions on the distribution\footnote{The authors additionally require that the random variable $\eta(X)$ on $[0, 1]$ admits bounded density.} of the pair $(X, Y) \in \bbR^d \times \{0, 1\}$ they proposed an algorithm whose rate of convergence is at most $\bigO(a_n^{ - {(1 + 1 \wedge \alpha)}/{2}})$.
    This rate is rather counter intuitive, since it suggests that if the constant $\alpha$ in the margin assumption is large it does not affect the rate of convergence.
    In contrast, here we show that for the proposed algorithm the rate of convergence is of order $\bigO(a_n^{ - {(1 + \alpha)}/{2}})$.
    That is, it strictly improves upon the results in~\cite{Yan_Koyejo_Zhong_Ravikumar18} whenever the constant $\alpha > 1$.
    However, it should be noted, that the authors of~\cite{Yan_Koyejo_Zhong_Ravikumar18} study a much more general family of the score functions and the sub-optimal rate can result from such a generality.
    \item We show that the constructed estimator is optimal in the minimax sense over the class of H\"older smooth regression functions.
    Let us mention that the optimality of the bound is expected, as in the classical work of~\cite{Audibert_Tsybakov07} the authors showed that the minimax risk in the standard binary classification settings is of order $a_n^{ - {(1 + \alpha)}/{2}}$, and it is achieved by a plug-in rule classifier.
    Clearly, it is hard to expect that the rate in a more difficult situation can be improved.
    Nevertheless, to the best of our knowledge, the minimax optimality in the context of binary classification with $\text{F}$-score have not been considered before.
\end{itemize}

The paper is organized as follows: in Section~\ref{sec:main_results} we present the semi-supervised classification algorithm; in Section~\ref{sec:upper} we establish an upper bound on the excess F-score under the margin assumption; in Section~\ref{sec:lower_bound} we introduce the class of distributions considered in this work and establish a minimax lower bound on the excess F-score.

%%%%%%%%%%%%%%%%%%%%%%%%%%%%%%%%%%%%%%%%%%%%%%%%%%%%%%%%%%%%%%%%%%%%%%%%%%%%%%%
\section{Main results}
\label{sec:main_results}
%%%%%%%%%%%%%%%%%%%%%%%%%%%%%%%%%%%%%%%%%%%%%%%%%%%%%%%%%%%%%%%%%%%%%%%%%%%%%%%
In this section we describe the proposed procedure $\hat g$ to estimate the Bayes optimal classifier $g^*$, this procedure is performed in two steps.
On the first step we estimate the regression function $\eta: \bbR^d \mapsto \{0, 1\}$ using the labeled data $\data_n$ and on the second step we estimate the optimal threshold $\theta^*$ based on the unlabeled data $\data_N$ and the estimator $\heta$ provided by the first step.
This procedure falls into the category of plug-in type classifiers, that is, we formally replace all the unknown quantities in the Bayes rule by its estimates.
That is, the classifier $\hat g$ is defined as
\begin{align*}
    \hat g(x) = \ind{\heta(x) > \htheta}\enspace,
\end{align*}
where $\heta$ is any estimator satisfying Assumption~\ref{ass:exponential_concentration_estimator} and $\htheta$ is the unique solution of
\begin{align}
    \label{eq:definition_threshold_estimator}
    \theta \frac{1}{N}\sum_{X_i \in \data_N} \heta(X_i) = \frac{1}{N}\sum_{X_i \in \data_N} (\heta(X_i) - \theta)_+\enspace.
\end{align}
In practice one can use a simple bisection algorithm~\cite[Algorithm 3.1]{Conte_Boor80} or its more sophisticated modifications (regula falsi or the secant method) to approximate $\htheta$ with any given precision.
For our theoretical analysis we assume that Equation~\eqref{eq:definition_threshold_estimator} is solved exactly.
However a simple modification of our arguments can handle the situation when the threshold $\htheta$ is known up to an additive factor $\epsilon_n = \bigO(a_n^{-1/2})$.

%%%%%%%%%%%%%%%%%%%%%%%%%%%%%%%%%%%%%%%%%%%%%%%%%%%%%%%%%%%%%%%%%%%%%%%%%%%%%%%
\subsection{Upper bound}
\label{sec:upper}
%%%%%%%%%%%%%%%%%%%%%%%%%%%%%%%%%%%%%%%%%%%%%%%%%%%%%%%%%%%%%%%%%%%%%%%%%%%%%%%

The main result of this subsection is an upper bound on excess score of the proposed procedure.
Here we provide two theorems, the first one gives an upper bound on the expected difference between the optimal threshold $\theta^*$ and its estimate $\htheta$.
The second one gives an upper bound on the excess $\text{F}$-score.
\begin{theorem}
    \label{thm:thresholds_upper_bound}
    If there exists an estimator $\heta$ of the regression function $\eta$ which satisfies Assumption~\ref{ass:exponential_concentration_estimator}, then there exists a constant $C > 0$ which depends on $C_0, C_1, C_2, p$ such that, the threshold $\htheta$ defined in Eq.~\eqref{eq:definition_threshold_estimator} satisfies
    \begin{align*}
        \Exp_{(\data_n, \data_N)}{|\theta^* - \htheta|} \leq C\left(a_n^{-1/2} + N^{-1/2}\right)\enspace.
    \end{align*}
\end{theorem}
\begin{theorem}
    \label{thm:upper_bound}
    If the distribution $\Prob$ of $(X, Y)$ satisfies the $\alpha$-margin assumption for some $C_0 > 0$ and $\alpha \geq 0$ and there exists an estimator $\heta$ of the regression function $\eta$ which satisfies Assumption~\ref{ass:exponential_concentration_estimator}, then there exists a constant $C > 0$ which depends on $\alpha, C_0, C_1, C_2, p$ such that
    \begin{align*}
        \Exp_{(\data_n, \data_N)}\excess_1(\hat g) \leq C\left(a_n^{-\tfrac{1 + \alpha}{2}} + N^{-\tfrac{1 + \alpha}{2}}\right)\enspace,
    \end{align*}
    where $\hat g(x) = \ind{\heta(x) > \htheta}$ with the threshold $\htheta$ defined in Equation~\eqref{eq:definition_threshold_estimator}.
\end{theorem}
Before proceeding to the proofs let us discuss the implications of these results.
First of all, there are two regimes in the bound of Theorems~\ref{thm:upper_bound}, the first one is $N \geq a_n$, in this regime, the dominant term is $a_n^{-{(1 + \alpha)}/{2}}$ which is the classical rate of convergence in the standard settings of binary classification with the $\alpha$-margin assumption.
The second regime is when $N < a_n$, then the dominating term of the bound is $N^{-{(1 + \alpha)}/{2}}$.
However, let us recall that one can always augment the second unlabeled dataset $\data_N$ by dividing $\data_n$ into two independent parts.
It implies that the second regime never occurs in our theoretical analysis of the excess score and the upper bound is actually independent of $N$.
Similar reasoning holds for the case of the optimal threshold estimation in Theorem~\ref{thm:thresholds_upper_bound}.
Once it is clear that the obtained upper bounds are actually independent of the size of the unlabeled dataset $\data_N$ it is interesting to notice that the dependence on $n$ is the same as in the standard case of the binary classification~\cite{Audibert_Tsybakov07}.
That is, similarly to the standard settings, the binary classification with F-score can achieve fast (faster than $1/\sqrt{n}$) and even super-fast (faster than $1/n$) rate depending on the interplay of $\alpha, \beta, d$.

Proofs of both theorems relies on the following lemma, provided in Appendix~\ref{app:difference_of_thresholds}, which relates the difference of the thresholds to the difference of the cumulative distribution function empirical of (CDF) $\eta$ and empirical CDF of $\heta$.
\begin{lemma}
    \label{lem:difference_of_thresholds}
    Let $\htheta \in [0, 1]$ be the threshold which satisfies Equation~\ref{eq:definition_threshold_estimator}, then
    \begin{align*}
        \abs{\htheta - \theta^*}\Prob(Y = 1) \leq \int_{0}^1 \abs{\Prob_X(\eta(X) \leq t) - \frac{1}{N}\sum_{X_i \in \data_N}\ind{\heta(X_i) \leq t}} dt\enspace.
    \end{align*}
\end{lemma}
This result is the main reason why our conclusions on the semi-supervised estimation is different from the ones in~\cite{Chzhen_Denis_Hebiri19,Singh_Nowak_Zhu09}.
For instance, in~\cite{Chzhen_Denis_Hebiri19} the authors also obtain a final decision rule by thresholding on some estimated level.
However, in the present work the difference between $\theta^*$ and $\htheta$ is controlled via $\ell_1$-norm of difference of CDF's, whereas in~\cite{Chzhen_Denis_Hebiri19} they control a similar quantity through Wassertstein infinity distance.

The complete proof of Theorems~\ref{thm:thresholds_upper_bound}  and~\ref{thm:upper_bound} can be found in Appendix~\ref{app:upper}, we only sketch the steps which are different from the analysis of~\cite{Audibert_Tsybakov07}.
Recall, that due to Lemma~\ref{lem:excess_score} we have the following bound for the excess score $\excess_1$
\begin{align*}
    \Exp_{(\data_n, \data_N)} \frac{\Exp\abs{\eta(X) - \theta^*}\ind{g^*(X) \neq \hat g (X)}}{\Prob(Y = 1) + \Prob(\hat g(X) = 1)}
    \leq \frac{1}{p}\Exp_{(\data_n, \data_N)}\Exp\abs{\eta(X) - \theta^*}\ind{g^*(X) \neq \hat g (X)}\enspace.
\end{align*}
First of all, notice that if for some $x \in \bbR^d$ the event $g^*(x) \neq \hat g (x)$ occurs, than we have
\begin{align*}
    \abs{\eta(x) - \theta^*} \leq \abs{\eta(x) - \heta(x)} + \absin{\theta^* - \htheta}\enspace,
\end{align*}
which further implies that at least one of the following inequalities hold for this $x \in \bbR^d$
\begin{align*}
    \abs{\eta(x) - \theta^*} &\leq 2\abs{\eta(x) - \heta(x)}\enspace,\\
    \abs{\eta(x) - \theta^*} &\leq 2\absin{\theta^* - \htheta}\enspace.
\end{align*}
Thus, we can upper bound the excess risk as
\begin{align*}
    \excess_1(\hat g) \leq  &\underbrace{\frac{1}{p}{\Exp\abs{\eta(X) - \theta^*}\ind{\abs{\eta(X) - \theta^*} \leq 2\abs{\eta(X) - \heta(X)}}}}_{T_1}\\
    &+
    \underbrace{\frac{1}{p}{\Exp\abs{\eta(X) - \theta^*}\ind{\abs{\eta(X) - \theta^*} \leq 2\abs{\theta^* - \htheta}}}}_{T_2}\enspace.
\end{align*}
The first term on the right hand side ($T_1$) of the inequality can be handled by the peeling technique used in~\cite[Lemma 3.1.]{Audibert_Tsybakov07},
which implies that, there exists a constant $C' = C'(p, \alpha, C_0, C_1, C_2) > 0$ such that
\begin{align*}
    \Exp_{(\data_n, \data_N)}T_1 \leq C'a_n^{-\tfrac{1 + \alpha}{2}}\enspace.
\end{align*}
Hence, it remains to upper bound the second term on the right hand side $(T_2)$ of the inequality.
Using Lemma~\ref{lem:difference_of_thresholds} we can upper bound $T_2$ as
\begin{align*}
    T_2 \leq \frac{1}{p}{\Exp\abs{\eta(X) - \theta^*}\ind{E}}\enspace,
\end{align*}
with $E = \ens{p\abs{\eta(X) - \theta^*} \leq 2\int_{0}^1 \abs{\Prob_X(\eta(X) \leq t) - \frac{1}{N}\sum_{X_i \in \data_N}\ind{\heta(X_i) \leq t}} dt}$.
Finally, we upper bound the indicator $\ind{E}$ by the indicators of two events $E^1$ and $E^2$ which are defined as
\begin{align*}
    E^1
    &=
    \ens{p\abs{\eta(X) - \theta^*} \leq 4\sup_{t \in [0, 1]} \abs{\Prob_X(\heta(X) \leq t) - \frac{1}{N}\sum_{X_i \in \data_N}\ind{\heta(X_i) \leq t}}}\enspace,\\
    E^2
    &=
    \ens{p\abs{\eta(X) - \theta^*} \leq 4\int_{0}^1 \abs{\Prob_X(\heta(X) \leq t) - \Prob_X(\eta(X) \leq t)} dt}\enspace.
\end{align*}
Thus, we have the following upper bound on $T_2$
\begin{align*}
    T_2
    \leq
    \underbrace{\frac{1}{p}{\Exp\abs{\eta(X) - \theta^*}\ind{E^1}}}_{T_2^1}
    +
    \underbrace{\frac{1}{p}{\Exp\abs{\eta(X) - \theta^*}\ind{E^2}}}_{T_2^2}\enspace,
\end{align*}
Notice that thanks to the Dvoretzky-Kiefer-Wolfowitz inequality~\cite{Dvoretzky_Kiefer_Wolfowitz56,Massart90} the term
\begin{align*}
    \sup_{t \in [0, 1]} \abs{\Prob_X(\heta(X) \leq t) - \frac{1}{N}\sum_{X_i \in \data_N}\ind{\heta(X_i) \leq t}}\enspace,
\end{align*}
conditionally on $\data_n$ admits an exponential concentration with the rate $N^{-1/2}$.
Hence, using the margin assumption, one can effortlessly show there exists a constant $C'' = C''(p, \alpha, C_0) > 0$ such that
 \begin{align*}
     \Exp_{(\data_n, \data_N)}T_2^1 \leq C''N^{- \tfrac{1 + \alpha}{2}}\enspace.
 \end{align*}
For the second term $T_2^2$ we proceed as follows
\begin{align*}
    &T_2^2 \leq \frac{4}{p^2}\int_{0}^1 \abs{\Prob_X(\heta(X) \leq t) - \Prob_X(\eta(X) \leq t)} dt \Prob({E^2})\enspace,
\end{align*}
thus, using the $\alpha$-margin assumption we get
\begin{align*}
    T_2^2 \leq\frac{C_04^{1 + \alpha}}{p^{2 + \alpha}}\left(\int_{0}^1 \abs{\Prob_X(\heta(X) \leq t) - \Prob_X(\eta(X) \leq t)} dt\right)^{1 + \alpha}\enspace,
\end{align*}
the integral on the right hand side of the bound corresponds to the $1$-Wasserstein distance on the real line, see for instance~\cite[Theorem 2.9]{Bobkov_Ledoux16} or~\cite{Vallender74} for the proof, and can be further upper bounded by the $L_1$ norm between $\heta$ and $\eta$, that is
\begin{align*}
    T_2^2 \leq\frac{C_04^{1 + \alpha}}{p^{2 + \alpha}}\left(\Exp_{\Prob_X}\abs{\eta(X) - \heta(X)}\right)^{1 + \alpha}\enspace.
\end{align*}
% In order to upper bound the expectation of $T_2^2$ we notice that the function $x \mapsto x^{1 + \alpha}$ is convex for all non-negative values of $\alpha$, then using Jensen's inequality we can write
% \begin{align*}
%     \Exp_{(\data_n, \data_N)}T_2^2 = \Exp_{\data_n}T_2^2 \leq \frac{C_04^{1 + \alpha}}{p^{2 + \alpha}}\Exp_{\Prob_X}\Exp_{\data_n}\abs{\eta(X) - \heta(X)}^{1 + \alpha}\enspace,
% \end{align*}
Since the estimator $\heta$ satisfies Assumption~\ref{ass:exponential_concentration_estimator}, one can show that there exists a constant $C''' = C'''(p, \alpha, C_0, C_1, C_2) > 0$ such that
\begin{align*}
    \Exp_{(\data_n, \data_N)}T_2^2 \leq C'''a_n^{-\tfrac{1 + \alpha}{2}}\enspace.
\end{align*}
Combination of all the inequalities yields the result of Theorem~\ref{thm:upper_bound}.
Notice that the same reasoning starting from Lemma~\ref{lem:difference_of_thresholds} implies the upper bound on the threshold estimation, that is, Theorem~\ref{thm:thresholds_upper_bound}.
    % There are several implications of the upper bound above.
    % On the one hand, notice that in the case $N \geq n$, the dominant term in the upper bound is of order $a_n^{-\tfrac{1 + \alpha}{2}}$, it suggests that it is 'harder' to estimate the regression function $\eta: \bbR^d \mapsto [0, 1]$, than the optimal threshold $\theta^*$.
    % On the other hand, the regime when $N < n$ is not interesting since one can always separate  the labeled dataset $\data_n$ into two independent datasets and use one them to estimate $\eta$ and the other one to estimate $\theta^*$.
    % That is, for any $N, n \in \bbN$ the minimax rate is at most $O\left(a_n^{-\tfrac{1 + \alpha}{2}}\right)$.

%%%%%%%%%%%%%%%%%%%%%%%%%%%%%%%%%%%%%%%%%%%%%%%%%%%%%%%%%%%%%%%%%%%%%%%%%%%%%%%
\subsection{Lower bound}
\label{sec:lower_bound}
%%%%%%%%%%%%%%%%%%%%%%%%%%%%%%%%%%%%%%%%%%%%%%%%%%%%%%%%%%%%%%%%%%%%%%%%%%%%%%%
In the beginning of the section we state the class of distribution $\class{P}_{\Sigma}$ of the random pair $(X, Y) \in \bbR^d \times \{0, 1\}$ considered in this work.
The first assumption is made on smoothness of the regression function $\eta: \bbR^d \mapsto [0, 1]$.
\begin{definition}[H\"older smoothness]
    Let $L > 0$ and $\beta > 0$. The class of function $\Sigma(\beta, L, \bbR^d)$ consists of all functions $h: \bbR^d \mapsto [0, 1]$ such that for all $x, x' \in \bbR^d$, we have
    \begin{align*}
        \abs{h(x) - h_x(x')} \leq L \norm{x - x'}^\beta_2\enspace,
    \end{align*}
    where $h_x(\cdot)$ is the Taylor polynomial of $h$ at point $x$ of degree $\floor{\beta}$.
\end{definition}
\begin{assumption}[$(\beta, L)$-H\"older regression function]
    \label{ass:holder}
    The distribution $\Prob$ of the pair $(X, Y) \in \bbR^d \times \{0, 1\}$ is such that $\eta \in \Sigma(\beta, L, \bbR^d)$ for some positive $\beta, L$.
\end{assumption}
Assumption~\ref{ass:holder} is usually not sufficient to guarantee the existence of an estimator $\heta$ satisfying Assumption~\ref{ass:exponential_concentration_estimator}:
extra assumptions are required on the marginal distribution $\Prob_X$ of the vector $X \in \bbR^d$.
\begin{definition}
    A Lebesgue measurable set $A \subset \bbR^d$ is said to be $(c_0, r_0)$-regular for some constants $c_0 > 0, r_0 > 0$ if for every $x \in A$ and every $r \in (0, r_0]$ we have
    \begin{align*}
        \lambda\left(A \cap \mathcal{B}(x, r)\right) \geq c_0\lambda\left(\mathcal{B}(x, r)\right)\enspace,
    \end{align*}
    where $\lambda$ is the Lebesgue measure and $\mathcal{B}(x, r)$ is the Euclidean ball of radius $r$ centered at $x$.
\end{definition}
\begin{assumption}[Strong density assumption]
    \label{ass:strong_density}
    We say that the marginal distribution $\Prob_X$ of the vector $X \in \bbR^d$ satisfies the strong density assumption if
    \begin{itemize}
        \item $\Prob_X$ is supported on a compact $(c_0, r_0)$-regular set $A \subset \bbR^d$,
        \item $\Prob_X$ admits a density $\mu$ \wrt to the Lebesgue measure uniformly lower- and upper-bounded by ${\mu_{\min}} > 0$ and $\mu_{\max} > 0$ respectively.
    \end{itemize}
\end{assumption}
If the regression function $\eta: \bbR^d \mapsto [0, 1]$ is $(\beta, L)$-H\"older and the marginal distribution satisfies the strong density assumption, one can state the following result due to~\cite{Audibert_Tsybakov07}.
\begin{theorem}[\cite{Audibert_Tsybakov07}]
    \label{thm:exponential_concentration}
    Let $\class{P}$ be a class of distributions on $\bbR^d \times \{0, 1\}$ such that the regression function $\eta \in \Sigma(\beta, L, \bbR^d)$ and the marginal distribution $\Prob_X$ satisfies the strong density assumption.
    Then, there exists an estimator $\heta$ of the regression function satisfying
    \begin{align*}
        \sup_{\Prob \in \class{P}} \Prob^{\otimes n}(\abs{\heta(x) - \eta(x)} \geq t) \leq C_1\exp\left(-C_2n^{\tfrac{2\beta}{2\beta + d}} t^2\right) \text{ a.s. } \Prob_X\enspace,
    \end{align*}
    for come constants $C_1, C_2$ depending on $\beta, d, L, c_0, r_0$.
\end{theorem}
Consider a class of distribution $\class{P}_{\Sigma}$ for which Assumptions~\ref{ass:margin_assumption},~\ref{ass:lower_proba_1},~\ref{ass:holder},~\ref{ass:strong_density} are satisfied, then Theorem~\ref{thm:exponential_concentration} and Theorems~\ref{thm:thresholds_upper_bound},~\ref{thm:upper_bound} imply the following corollary.
\begin{corollary}
    There exist constants $C, B > 0$ which depend only on $\alpha, p, d, C_0, C_1, C_2$ such that for any $n > 1, N > 1$ we have
    \begin{align}
        \inf_{\hat g}\sup_{\Prob \in \class{P}_{\Sigma}} \Exp_{(\data_n, \data_N)}\excess_1(\hat g)
        &\leq C n^{-\frac{(1 + \alpha)\beta}{2\beta + d}}\enspace,\\
        \inf_{\hat \theta}\sup_{\Prob \in \class{P}_{\Sigma}} \Exp_{(\data_n, \data_N)}\abs{\theta^* - \htheta}
        &\leq B n^{-\frac{\beta}{2\beta + d}}\enspace.
    \end{align}
    where the infima are taken over all estimators $\hat g$ and $\htheta$ respectively.
\end{corollary}
The next theorem states that the upper bounds of the previous corollary are optimal up to a constant multiplicative factor.
\begin{theorem}
    If $\alpha \beta \leq d$, there exists constants $c> 0$ such that for any $n > 1, N > 1$ we have the following lower-bound on the minimax risk
    \begin{align}
        \inf_{\hat g}\sup_{\Prob \in \class{P}_{\Sigma}} \Exp_{(\data_n, \data_N)}\excess_1(\hat g) &\geq c n^{-\frac{(1 + \alpha)\beta}{2\beta + d}}\enspace,
        % \inf_{\hat \theta}\sup_{\Prob \in \class{P}_{\Sigma}} \Exp_{(\data_n, \data_N)}\abs{\theta^* - \htheta}
        % &\geq b n^{-\frac{\beta}{2\beta + d}}\enspace,
    \end{align}
    where the infimum is taken over all estimators $\hat g$.
\end{theorem}
The proof of the lower bound can be found in Appendix~\ref{app:lower}, it follows standard information-theoretic arguments using reduction of the minimax risk to a Bayes risk.
The construction of the distributions is inspired by both~\cite{Rigollet_Vert09} and~\cite{Audibert_Tsybakov07}, and the actual proof relies on~\cite[Lemma 5.1.]{Audibert04}, which is based on the Assouad's lemma, see for instance~\cite[Lemma 2.12]{Tsybakov09}.
%%%%%%%%%%%%%%%%%%%%%%%%%%%%%%%%%%%%%%%%%%%%%%%%%%%%%%%%%%%%%%%%%%%%%%%%%%%%%%%
\section{Conclusion}
%%%%%%%%%%%%%%%%%%%%%%%%%%%%%%%%%%%%%%%%%%%%%%%%%%%%%%%%%%%%%%%%%%%%%%%%%%%%%%%
In this work we proposed a semi-supervised plug-in type algorithm for the problem of binary classification with $\text{F}$-score.
The proposed algorithm can leverage an unlabeled dataset for the estimation of the optimal threshold.
Under the margin assumption it is shown that the proposed algorithm is optimal in the minimax sense and can achieve fast rates of convergence.
Further development of the binary classification with F-score will be devoted to empirical risk minimization rules.
% \acks{We thank Mo7 and Joe for being cool.}
\begin{acknowledgements}
    This work was partially supported by ``Labex B\'ezout'' of Universit\'e Paris-Est.
    Besides, we would like to thank Joseph Salmon and Mohamed Hebiri for their thoughtful remarks.
\end{acknowledgements}

% \bibliography{references_all}

\appendix

\section{Bayes classifier and Lemma~\ref{lem:excess_score}}
\label{app:bayes_excess}

For the rest of this section the parameter $b > 0$ is assumed to be fixed and known.
Let us first recall the definition of the $\text{F}_b$-score
\begin{align*}
    F_b(g) = (1 + b^2)\frac{\Prob(Y = 1, g(X) = 1)}{b^2\Prob(Y = 1) + \Prob(g(X) = 1)}\enspace,
\end{align*}
and an optimal classifier is defined as
\begin{align*}
    g^* \in \argmax_{g \in \class{G}}F_b(g)\enspace.
\end{align*}
In this section we would like to show that a classifier defined for all $x \in \bbR^d$ as
\begin{align*}
    g_*(x) = \ind{\eta(x) \geq \theta^*}\enspace,
\end{align*}
with $\theta^*$ being a root of
\begin{align*}
   \theta \mapsto b^2 \Prob(Y = 1) \theta - \Exp(\eta(X) - \theta)_+\enspace.
\end{align*}
Let us first show that $\theta^*$ is well-defined, that is, it exists and is unique for every distribution with $\Prob(Y = 1) \neq 0$.
Hence, we would like to study solutions of the following equation
\begin{align*}
    b^2 \Prob(Y = 1) \theta =  \Exp(\eta(X) - \theta)_+\enspace.
\end{align*}
Clearly, the mapping $\theta \mapsto b^2 \Prob(Y = 1) \theta$ is continuous and strictly increasing on $[0, 1]$ and the mapping $\theta \mapsto \Exp(\eta(X) - \theta)_+$ is non-increasing on $[0, 1]$.
Thus, it is sufficient to demonstrate that the mapping $\theta \mapsto \Exp(\eta(X) - \theta)_+$ is continuous, indeed, let $\theta, \theta' \in [0, 1]$, then, due to the Lipschitz continuity of $(\cdot)_+$ we can write
\begin{align*}
   \abs{ \Exp(\eta(X) - \theta)_+ - \Exp(\eta(X) - \theta')_+} \leq \Exp\abs{(\eta(X) - \theta)_+ - (\eta(X) - \theta')_+} \leq \abs{\theta - \theta'}\enspace.
\end{align*}
This implies that the mapping $\theta \mapsto \Exp(\eta(X) - \theta)_+$ is a contraction and thus is continuous.
Hence, the threshold $\theta^*$ is well-defined, that is, it exists and is unique.
Consequently, the classifier $x \mapsto \ind{\eta(x) \geq \theta^*}$ is well-defined.

Now, we are interested in the value $F_b(g_*)$, we can write
\begin{align*}
    F_b(g_*)
    &=
    (1 + b^2)\frac{\Prob(Y = 1, g_*(X) = 1)}{b^2\Prob(Y = 1) + \Prob(g_*(X) = 1)}\\
    &=
    (1 + b^2)\frac{\Exp[\eta(X)\ind{\eta(X) \geq \theta^*}]}{b^2\Prob(Y = 1) + \Prob(\eta(X) \geq \theta^*)}\\
    &=
    (1 + b^2)\frac{\Exp[(\eta(X) - \theta^*)\ind{\eta(X) \geq \theta^*}] + \theta^*\Exp\ind{\eta(X) \geq \theta^*}}{b^2\Prob(Y = 1) + \Prob(\eta(X) \geq \theta^*)}\\
    &=
    (1 + b^2)\frac{\Exp(\eta(X) - \theta^*)_+ + \theta^*\Prob(\eta(X) \geq \theta^*)}{b^2\Prob(Y = 1) + \Prob(\eta(X) \geq \theta^*)}\enspace,
\end{align*}
using the definition of $\theta^*$ we continue as
\begin{align*}
    F_b(g_*)
    &=
    (1 + b^2)\frac{\Exp(\eta(X) - \theta^*)_+ + \theta^*\Prob(\eta(X) \geq \theta^*)}{b^2\Prob(Y = 1) + \Prob(\eta(X) \geq \theta^*)}\\
    &=
    (1 + b^2)\frac{\theta^*b^2\Prob(Y = 1) + \theta^*\Prob(\eta(X) \geq \theta^*)}{b^2\Prob(Y = 1) + \Prob(\eta(X) \geq \theta^*)} = (1 + b^2)\theta^*\enspace.
\end{align*}
To conclude the optimality of $g_*$ we prove Lemma~\ref{lem:excess_score}.
\begin{proof}
    Fix an arbitrary measurable function $g: \bbR^d \mapsto \{0, 1\}$, then by the definition of the excess score we have
    \begin{align*}
        \excess_b(g)
        &\coloneqq
        \frac{\Prob(Y = 1, g^*(X) = 1)}{b^2\Prob(Y = 1) + \Prob(g^*(X) = 1)} - \frac{\Prob(Y = 1, g(X) = 1)}{b^2\Prob(Y = 1) + \Prob(g(X) = 1)}\\
        &=
        \frac{\Exp\eta(X)\ind{\eta(X) > \theta^*}}{b^2\Prob(Y = 1) + \Prob(g^*(X) = 1)} - \frac{\Exp\eta(X)\ind{g(X) = 1}}{b^2\Prob(Y = 1) + \Prob(g(X) = 1)}\\
        &=
        \frac{\Exp\eta(X)\ind{\eta(X) > \theta^*} - \Exp\eta(X)\ind{g(X) = 1}}{b^2\Prob(Y = 1) + \Prob(g^*(X) = 1)}\\
        &\phantom{=}
        +
        \frac{\Exp\eta(X)\ind{g(X) = 1}}{b^2\Prob(Y = 1) + \Prob(g(X) = 1)}\left(\frac{\Prob(g(X) = 1) - \Prob(g^*(X) = 1)}{b^2\Prob(Y = 1) + \Prob(g^*(X) = 1)}\right)\\
        &=
        \frac{\Exp(\eta(X) - \theta^*)(\ind{\eta(X) > \theta^*} - \ind{g(X) = 1}) + \theta^*(\Prob(g^*(X) = 1) - \Prob(g(X) = 1))}{b^2\Prob(Y = 1) + \Prob(g^*(X) = 1)}\\
        &\phantom{=}
        +
        F_b(g)\left(\frac{\Prob(g(X) = 1) - \Prob(g^*(X) = 1)}{b^2\Prob(Y = 1) + \Prob(g^*(X) = 1)}\right)\\
        &=
        \frac{\Exp\abs{\eta(X) - \theta^*}\ind{g^*(X) \neq g(X)}}{b^2\Prob(Y = 1) + \Prob(g^*(X) = 1)} + (\theta^* - F_b(g))\frac{\Prob(g^*(X) = 1) - \Prob(g(X) = 1)}{b^2\Prob(Y = 1) + \Prob(g^*(X) = 1)}\enspace.
    \end{align*}
    Using Theorem~\ref{thm:bayes_optimal} we know that $\theta^* = F_b(g^*)$ and therefore
    \begin{align*}
        \excess(g) = \frac{\Exp\abs{\eta(X) - \theta^*}\ind{g^*(X) \neq g(X)}}{b^2\Prob(Y = 1) + \Prob(g^*(X) = 1)} + \excess(g)\frac{\Prob(g^*(X) = 1) - \Prob(g(X) = 1)}{b^2\Prob(Y = 1) + \Prob(g^*(X) = 1)}\enspace.
    \end{align*}
    We conclude by solving the previous equality for $\excess(g)$.
    Thus, $g_*$ is a Bayes optimal classifier and hence can be denoted by $g^*$.
\end{proof}

\section{Proof of Lemma~\ref{lem:difference_of_thresholds}}
\label{app:difference_of_thresholds}

\begin{proof}
    To prove this lemma, it is convenient to rewrite Equation~\ref{eq:definition_threshold_estimator} in terms of CDF.
    Let $\mu$ be an arbitrary probability measure on $\bbR^d$ and $p: \bbR^d \mapsto [0, 1]$ be any measurable function, then using Fubini's theorem we can write
    \begin{align*}
        \int p(x) d\mu(x)
        &=
        \int \int_{0}^1 \ind{p(x) > t} dt d\mu(x)\\
        &=
        \int_{0}^1 \mu(p(X) > t) dt\enspace,
    \end{align*}
    and for any $\theta \in [0, 1]$, since $(p(x) - \theta)_+ \in [0, 1]$ we have
    \begin{align*}
        \int (p(x) - \theta)_+ d\mu(x)
        &=
        \int \int_{0}^1 \ind{p(x) - \theta > t} dt d\mu(x)\\
        &=
        \int \int_{\theta}^{1 + \theta} \ind{p(x) > t} dt d\mu(x)\\
        &=
        \int \int_{\theta}^{1} \ind{p(x) > t} dt d\mu(x)\\
        &=
        \int_{\theta}^{1} \mu(p(X) > t) dt\enspace.
    \end{align*}
    Let us denote by $\Prob_{X, N} = \tfrac{1}{N}\sum_{X_i \in \data_N} \delta_{X_i}$ the empirical measure of the unlabeled dataset $\data_N$.
    Using these equalities, the thresholds $\theta^*, \htheta \in [0, 1]$ satisfy
    \begin{align*}
        \htheta = \frac{\int_{\htheta}^1 \Prob_{X, N}(\heta(X) > t) dt}{\int_{0}^1 \Prob_{X, N}(\heta(X) > t) dt},\quad \theta^* = \frac{\int_{\theta^*}^1 \Prob_{X}(\eta(X) > t) dt}{\int_{0}^1 \Prob_{X}(\eta(X) > t) dt}\enspace.
    \end{align*}
    Now, we are in position to bound the difference $|\htheta - \theta^*|$, first assume that $\theta^* \geq \htheta$, then
    \begin{align*}
        \theta^* - \htheta
        &=
        \frac{\int_{\theta^*}^1 \Prob_{X}(\eta(X) > t) dt}{\int_{0}^1 \Prob_{X}(\eta(X) > t) dt} -  \frac{\int_{\htheta}^1 \Prob_{X, N}(\heta(X) > t) dt}{\int_{0}^1 \Prob_{X, N}(\heta(X) > t) dt}\\
        &\leq
        \frac{\int_{\htheta}^1 \Prob_{X}(\eta(X) > t) dt}{\int_{0}^1 \Prob_{X}(\eta(X) > t) dt} -  \frac{\int_{\htheta}^1 \Prob_{X, N}(\heta(X) > t) dt}{\int_{0}^1 \Prob_{X, N}(\heta(X) > t) dt}\\
        &=
        \frac{\int_{\htheta}^1 (\Prob_{X}(\eta(X) > t) - \Prob_{X, N}(\heta(X) > t)) dt}{\int_{0}^1 \Prob_{X}(\eta(X) > t) dt}\\
        &\phantom{=}
        +
        \frac{\int_{\htheta}^1 \Prob_{X, N}(\heta(X) > t) dt}{\int_{0}^1 \Prob_{X, N}(\heta(X) > t) dt}\frac{\int_{0}^1 (\Prob_{X, N}(\eta(X) > t) - \Prob_{X}(\heta(X) > t)) dt}{\int_{0}^1 \Prob_{X}(\eta(X) > t) dt}\\
        &=
        \frac{\int_{\htheta}^1 (\Prob_{X}(\eta(X) > t) - \Prob_{X, N}(\heta(X) > t)) dt - \htheta \int_{0}^1 (\Prob_{X}(\heta(X) > t) - \Prob_{X, N}(\eta(X) > t)) dt}{\int_{0}^1 \Prob_{X}(\eta(X) > t) dt}\\
        &\leq
        \frac{1}{\Prob(Y  = 1)}\int_{0}^1\abs{\Prob_{X}(\eta(X) > t) - \Prob_{X, N}(\heta(X) > t)} dt\enspace.
    \end{align*}
    Further, if $\htheta > \theta^*$ we can write
    \begin{align*}
        \htheta - \theta^*
        &=
        \frac{\int_{\htheta}^1 \Prob_{X, N}(\heta(X) > t) dt}{\int_{0}^1 \Prob_{X, N}(\heta(X) > t) dt} - \frac{\int_{\theta^*}^1 \Prob_{X}(\eta(X) > t) dt}{\int_{0}^1 \Prob_{X}(\eta(X) > t) dt}\\
        &\leq
        \frac{\int_{\theta^*}^1 \Prob_{X, N}(\heta(X) > t) dt}{\int_{0}^1 \Prob_{X, N}(\heta(X) > t) dt} - \frac{\int_{\theta^*}^1 \Prob_{X}(\eta(X) > t) dt}{\int_{0}^1 \Prob_{X}(\eta(X) > t) dt}\\
        &=
        \frac{\int_{\theta^*}^1 (\Prob_{X, N}(\heta(X) > t) - \Prob_{X}(\eta(X) > t)) dt}{\int_{0}^1 \Prob_{X}(\eta(X) > t) dt}\\
        &\phantom{=}
        +
        \frac{\int_{\theta^*}^1 \Prob_{X, N}(\heta(X) > t) dt}{\int_{0}^1 \Prob_{X, N}(\heta(X) > t) dt}\frac{\int_{0}^1 (\Prob_{X}(\eta(X) > t) - \Prob_{X, N}(\heta(X) > t)) dt}{\int_{0}^1 \Prob_{X}(\eta(X) > t) dt}\\
        &\leq
        \frac{1}{\Prob(Y  = 1)}\int_{0}^1\abs{\Prob_{X}(\eta(X) > t) - \Prob_{X, N}(\heta(X) > t)} dt\enspace,
    \end{align*}
    where the last inequality follows the same lines as for the case $\htheta \leq \theta^*$.
\end{proof}

\section{Proof of the upper bound}
\label{app:upper}
Let $\hat\eta$ be an estimator of the regression function based on the labeled dataset $\data_n$ which satisfies Assumption~\ref{ass:exponential_concentration_estimator}.
Recall, that the estimator $\hat g$ is defined for every $x \in \bbR^d$ as
\begin{align*}
    \hat g(x) = \ind{\hat\eta(x) > \hat\theta}\enspace,
\end{align*}
with $\hat\theta$ being the unique solution of Eq.~\eqref{eq:definition_threshold_estimator}.
Unless stated otherwise, we work conditionally on $(\data_n, \data_N)$.
Using Lemma~\ref{lem:excess_score} we can express the excess score of $\hat g$ as
\begin{align*}
    \excess_1(\hat g)
    &=
    \frac{\Exp\abs{\eta(X) - \theta^*}\ind{g^*(X) \neq \hat g(X)}}{\Prob(Y = 1) + \Prob(\hat g(X) = 1)}
    \leq
    \frac{1}{p}\Exp\abs{\eta(X) - \theta^*}\ind{g^*(X) \neq \hat g(X)}\enspace.
\end{align*}
Clearly, on the event $\ens{g^*(X) \neq \hat g(X)}$ it holds that $\ens{\abs{\eta(X) - \theta^*} \leq \abs{\hat\eta(X) - \eta(X)} + \abs{\hat\theta - \theta^*}}$, thus
\begin{align*}
    \excess_1(\hat g)
    &\leq
    \frac{1}{p}\Exp\abs{\eta(X) - \theta^*}\ind{\abs{\eta(X) - \theta^*} \leq \abs{\hat\eta(X) - \eta(X)} + \abs{\hat\theta - \theta^*}}\\
    &\leq
    \frac{1}{p}\Exp\abs{\eta(X) - \theta^*}\ind{\abs{\eta(X) - \theta^*} \leq 2\abs{\hat\theta - \theta^*}}
    +
    \frac{1}{p}\Exp\abs{\eta(X) - \theta^*}\ind{\abs{\eta(X) - \theta^*} \leq 2\abs{\hat\eta(X) - \eta(X)}}\enspace.
\end{align*}
Using, Lemma~\ref{lem:difference_of_thresholds} the excess risk can be further upper bounded as
\begin{align*}
    \excess_1(\hat g)
    \leq
    &\frac{1}{p}\Exp\abs{\eta(X) - \theta^*}\ind{\abs{\eta(X) - \theta^*} \leq 2\abs{\hat\eta(X) - \eta(X)}}\\
    &+
    \frac{1}{p}\Exp\abs{\eta(X) - \theta^*}\ind{\abs{\eta(X) - \theta^*} \leq \frac{2}{p}\int_{0}^1 \abs{\Prob_X(\eta(X) \leq t) - \frac{1}{N}\sum_{X_i \in \data_N}\ind{\heta(X_i) \leq t}} dt}\\
    \leq
    &\frac{1}{p}\Exp\abs{\eta(X) - \theta^*}\ind{\abs{\eta(X) - \theta^*} \leq 2\abs{\hat\eta(X) - \eta(X)}}\\
    &+
    \frac{1}{p}\Exp\abs{\eta(X) - \theta^*}\ind{\abs{\eta(X) - \theta^*} \leq \frac{2}{p}\int_{0}^1 \abs{\Prob_X(\eta(X) \leq t) - \Prob_X(\hat\eta(X) \leq t)} dt}\\
    &+
    \frac{1}{p}\Exp\abs{\eta(X) - \theta^*}\ind{\abs{\eta(X) - \theta^*} \leq \frac{2}{p}\int_{0}^1 \abs{\frac{1}{N}\sum_{X_i \in \data_N}\ind{\heta(X_i) \leq t} - \Prob_X(\hat\eta(X) \leq t)} dt}\enspace.
\end{align*}
Notice that $\int_{0}^1 \abs{\Prob_X(\eta(X) \leq t) - \Prob_X(\hat\eta(X) \leq t)} dt = \norm{F_{\eta} - F_{\hat\eta}}_1$, with $F_\eta, F_{\hat \eta}$ being the cumulative distribution functions of $\eta, \hat\eta$ respectively, corresponds to the 1-Wasserstein distance, see~\cite{Bobkov_Ledoux16} for an in-depth discussion.
Therefore, we have
\begin{align*}
    \int_{0}^1 \abs{\Prob_X(\eta(X) \leq t) - \Prob_X(\hat\eta(X) \leq t)} dt \leq \Exp_{X \sim \Prob_X}\abs{\eta(X) - \hat\eta(X)} \eqdef\norm{\eta - \hat\eta}_1\enspace,
\end{align*}
 and introducing notation $\hat\Prob_X \eqdef \tfrac{1}{N}\sum_{X_i \in \data_N}\delta_{X_i}$ for the empirical measure of the feature vector $X$ we can write
\begin{align*}
    \excess_1(\hat g)
    \leq
    &\frac{1}{p}\Exp\abs{\eta(X) - \theta^*}\ind{\abs{\eta(X) - \theta^*} \leq 2\abs{\hat\eta(X) - \eta(X)}}\\
    &+
    \frac{1}{p}\Exp\abs{\eta(X) - \theta^*}\ind{\abs{\eta(X) - \theta^*} \leq \frac{2}{p}\norm{\eta - \hat\eta}_1}\\
    &+
    \frac{1}{p}\Exp\abs{\eta(X) - \theta^*}\ind{\abs{\eta(X) - \theta^*} \leq \frac{2}{p}\sup_{t \in [0, 1]} \abs{\hat\Prob_X(\hat\eta(X) \leq t) - \Prob_X(\hat\eta(X) \leq t)}}\enspace.
\end{align*}
Finally, using the margin Assumption~\ref{ass:margin_assumption} we can write
\begin{align*}
    \excess_1(\hat g)
    \leq
    &\frac{1}{p}\Exp\abs{\eta(X) - \theta^*}\ind{\abs{\eta(X) - \theta^*} \leq 2\abs{\hat\eta(X) - \eta(X)}}\\
    &+
    \frac{1}{p}\Exp\abs{\eta(X) - \theta^*}\ind{\abs{\eta(X) - \theta^*} \leq \frac{2}{p}\norm{\eta - \hat\eta}_1}\\
    &+
    \frac{1}{p}\Exp\abs{\eta(X) - \theta^*}\ind{\abs{\eta(X) - \theta^*} \leq \frac{2}{p}\sup_{t \in [0, 1]} \abs{\hat\Prob_X(\hat\eta(X) \leq t) - \Prob_X(\hat\eta(X) \leq t)}}\\
    \leq
    &\frac{1}{p}\Exp\abs{\eta(X) - \theta^*}\ind{\abs{\eta(X) - \theta^*} \leq 2\abs{\hat\eta(X) - \eta(X)}}\\
    &+
    \frac{2}{p^2}\norm{\eta - \hat\eta}_1\Prob\parent{\abs{\eta(X) - \theta^*} \leq \frac{2}{p}\norm{\eta - \hat\eta}_1}\\
    &+
    \frac{1}{p}\Exp\abs{\eta(X) - \theta^*}\ind{\abs{\eta(X) - \theta^*} \leq \frac{2}{p}\sup_{t \in [0, 1]} \abs{\hat\Prob_X(\hat\eta(X) \leq t) - \Prob_X(\hat\eta(X) \leq t)}}\\
    \leq
    &\frac{1}{p}\Exp\abs{\eta(X) - \theta^*}\ind{\abs{\eta(X) - \theta^*} \leq 2\abs{\hat\eta(X) - \eta(X)}}
    +
    \frac{2^{\alpha + 1}c_0}{p^{2 + \alpha}}\norm{\eta - \hat\eta}^{1 + \alpha}_1\\
    &+
    \frac{1}{p}\Exp\abs{\eta(X) - \theta^*}\ind{\abs{\eta(X) - \theta^*} \leq \frac{2}{p}\sup_{t \in [0, 1]} \abs{\hat\Prob_X(\hat\eta(X) \leq t) - \Prob_X(\hat\eta(X) \leq t)}}\enspace.
\end{align*}
Taking expectation from the both sides with respect to the distribution of $(\data_n, \data_N)$ we follow~\cite[Lemma 3.1]{Audibert_Tsybakov07} to bound the first term on the right hand side.
This peeling argument became classical in the literature and thus is omitted here.
Moreover, using Assumption~\ref{ass:exponential_concentration_estimator} the second term can be bounded with the same rate as the first term.
These arguments would imply that there exists $C \geq 0$ such that for all $n, N \geq 1$ it holds that
\begin{align*}
    \Exp_{(\data_n, \data_N)}\excess_1(\hat g)
    &\leq
    C a_n^{-\frac{1 + \alpha}{2}}\\
    &\phantom{=}+
    \frac{1}{p}\Exp_{(\data_n, \data_N)}\Exp\abs{\eta(X) - \theta^*}\ind{\abs{\eta(X) - \theta^*} \leq \frac{2}{p}\sup_{t \in [0, 1]} \abs{\hat\Prob_X(\hat\eta(X) \leq t) - \Prob_X(\hat\eta(X) \leq t)}}\\
    &\leq
    C a_n^{-\frac{1 + \alpha}{2}}
    +
    \frac{2^{\alpha + 1}c_0}{p^{2 + \alpha}}\Exp_{(\data_n, \data_N)}\parent{\sup_{t \in [0, 1]} \abs{\hat\Prob_X(\hat\eta(X) \leq t) - \Prob_X(\hat\eta(X) \leq t)}}^{1 + \alpha}
\end{align*}
It remains to upper bound the second term in the bound above, to this end we recall the classical Dvoretzky-Kiefer-Wolfowitz inequality~\cite{Massart90}
\begin{lemma}[Dvoretzky-Kiefer-Wolfowitz inequality]
    Given $N \geq 0$, let $Z_1, \ldots, Z_N$ be \iid real-valued random variables with cumulative distribution function $F_Z$, denote by $\hat{F}_Z$ the cumulative distribution function with respect to the empirical measure, that is, with respect to $\frac{1}{N}\sum_{i = 1}^N\delta_{Z_i}$, then for every $t > 0$ we have
    \begin{align*}
        \Prob\parent{\sup_{z \in \bbR}\abs{\hat F_Z(z) - F_Z(z)} \geq t} \leq 2\exp\parent{-2Nt^2}\enspace.
    \end{align*}
\end{lemma}
Let us apply this lemma to $Z_i \eqdef \hat\eta(X_i)$, conditionally on $\data_n$ these random variables are \iid real-valued, thus for all $t > 0$
\begin{align*}
    \Prob\parent{\sup_{t \in [0, 1]} \abs{\hat\Prob_X(\hat\eta(X) \leq t) - \Prob_X(\hat\eta(X) \leq t} \geq t \Big\lvert \data_n} \leq  2\exp\parent{-2Nt^2}, \quad \text{a.s. } \data_n\enspace.
\end{align*}
Finally, to conclude the upper bound we apply this exponential concentration to upper bound the expectation as
\begin{align*}
    \Exp_{\data_n}\Exp_{\data_N}\left[\parent{\Delta_{(\data_N, \data_n)}}^{1 + \alpha} \Big\lvert \data_n\right]
    &=
    \Exp_{\data_n}\int_{0}^{\infty}\Prob\parent{{\Delta_{(\data_N, \data_n)}} \geq t^{\frac{1}{1 + \alpha}} \Big\lvert \data_n} dt\\
    &\leq
    \int_{0}^{\infty}2\exp\parent{-2Nt^{\frac{2}{1 + \alpha}}} dt\\
    &=
    N^{-\frac{1 + \alpha}{2}}2\int_{0}^{\infty}\exp\parent{-2t^{\frac{2}{1 + \alpha}}} dt\\
    &\leq
    CN^{-\frac{1 + \alpha}{2}}\enspace,
\end{align*}
where we used the shortcut $\Delta_{(\data_N, \data_n)}$ for the desired empirical process.
Combining all the bounds we conclude.
 \section{Proof of the lower bound}
\label{app:lower}

\begin{proof}
    The proof is similar to the one used in~\cite{Audibert_Tsybakov07} and in~\cite{Rigollet_Vert09} and is based on Assouad lemma.
    Similarly, we define the regular grid on $\bbR^d$ as
    \begin{align*}
        G_q \coloneqq \left\{\left(\frac{2k_1 + 1}{2q}, \ldots, \frac{2k_d + 1}{2q}\right)^\top \,:\, k_i \in \{0, \ldots, q - 1\}, i = 1, \ldots, d\right\}\enspace,
    \end{align*}
    and denote by $n_q(x) \in G_q$ as the closest point to of the grid $G_q$ to the point $x \in \bbR^d$.
    Such a grid defines a partition of the unit cube $[0, 1]^d \subset \bbR^d$ denoted by $\class{X}'_1, \ldots, \class{X}'_{q^d}$.
    Besides, denote by $\class{X}'_{-j} \coloneqq \{x \in \bbR^d\, : \, -x \in \class{X}_j'\}$ for all $j = 1, \ldots, q^d$.
    For a fixed integer $m \leq q^d$ and for any $j \in \{1, \ldots, m\}$ define $\class{X}_i \coloneqq \class{X}_i'$, $\class{X}_{-i} \coloneqq \class{X}_{-i}'$.
    For every $\sigma \in \{-1, 1\}^m$ we define a regression function $\eta_\sigma$ as
    \begin{align*}
        \eta_\sigma(x) = \begin{cases}
                            \frac{1}{4} + \sigma_j \varphi(x), &\text{ if } x \in \class{X}_i\\
                            \frac{1}{4} - \sigma_j \varphi(x), &\text{ if } x \in \class{X}_{-i}\\
                            \frac{1}{4}, &\text{ if } x \in \class{B}(0, \sqrt{d}) \setminus \left(\cup_{i = -m, i \neq 0}^m \class{X}_i\right)\\
                            \tau, &\text{ if } x \in \bbR^d \setminus \class{B}(0, \sqrt{d} + \rho)\\
                            \xi(x), &\text{ if } x \in \class{B}(0, \sqrt{d} + \rho) \setminus \class{B}(0, \sqrt{d})
                         \end{cases}\enspace,
    \end{align*}
    where $\rho, \varphi, \xi, \tau$ are to be specified and $\class{B}(0, \sqrt{d} + \rho), \class{B}(0, \sqrt{d})$ are Euclidean balls of radius $\sqrt{d} + \rho$ and $\sqrt{d}$ respectively.
    The definition of the function $\varphi$ is exactly the same as in~\cite{Audibert_Tsybakov07}.
    That is, $\varphi \coloneqq C_\varphi q^{-\beta} u(q\norm{x - n_q(x)}_2)$ with some non-increasing infinitely differentiable function such that $u(x) = 1$ for $x \in [0, 1/4]$ and $u(x) = 0$ for $x \geq 1/2$.
    The function $\xi$ is defined as $\xi(x) = (\tau - 1/4) v([\norm{x}_2 - \sqrt{d}] / \rho) + 1/4$, where $v$ is non-decreasing infinitely differentiable function such that $v(x) = 0$ for $x \leq 0$ and $v(x) = 1$ for $x \geq 1$.
    The constant $\rho$ is chosen big enough to ensure that $|\xi(x) - \xi_x(x')| \leq L \norm{x - x'}_2^{\beta}$ for any $x, x' \in \bbR^d$.

    For any $\sigma \in \{-1, 1\}^m$ we construct a marginal distribution $P_X$ which is independent of $\sigma$ and has a density $\mu$ \wrt to the Lebesgue measure on $\bbR^d$.
    Fix some $0 < w \leq m^{-1}$ and set $A_0$ a Euclidean ball in $\bbR^d$ that has an empty intersection with $\class{B}(0, \sqrt{d} + \rho)$ and whose Lebesgue measure is $\lambda(A_0) = 1 - mq^{-d}$.
    The density $\mu$ is constructed as
    \begin{itemize}
        \item $\mu(x) = \frac{w}{\lambda(\class{B}(0, (4q)^{-1}))}$ for every $z \in G_q$ and every $x \in \class{B}(z, (4q)^{-1}))$ or $x \in \class{B}(-z, (4q)^{-1}))$,
        \item $\mu(x) = \frac{1 - 2mw}{\lambda(A_0)}$ for every $x \in A_0$,
        \item $\mu(x) = 0$ for every other $x \in \bbR^d$.
    \end{itemize}
    To complete the construction it remain to specify the value of $\tau \in [0, 1]$.
    The idea here is to force the optimal threshold $\theta^*$ to be equal to some predefined constant using the additional degree of freedom provided by the parameter $\tau$.
    Importantly, this optimal threshold should not depend on the binary vector $\sigma \in \{-1, 1\}^m$.
    To achieve this recall that we set $\theta^* = 1/4$ and show that there exists an appropriate choice of $\tau$.
    First, recall that the optimal threshold $\theta^*$ satisfies
    \begin{align*}
        \theta^* \Exp \eta(X) = \Exp(\eta(X) - \theta^*)\enspace.
    \end{align*}
    Define $b' = \int_{\class{X}_1} \varphi(x) \mu(x) dx / \int_{\class{X}_1} \mu(x) dx$ and put $\theta^* = 1/4$, notice that the left hand side of the last equality for every $\sigma \in \{-1, 1\}^m$ is given by
    \begin{align*}
        \Exp_\mu \eta_\sigma(X) &= \int_{\bbR^d} \eta(x) d\mu(x)\\
                                &=
                                \sum_{j = 1}^m\int_{\class{X}_j} (1/4 + \sigma_j\xi(x)) d\mu(x)
                                +
                                \sum_{j = 1}^m\int_{\class{X}_{-j}} (1/4 - \sigma_j\xi(x)) d\mu(x)
                                +\int_{A_0}\tau d\mu(x)\\
                                & = \frac{mw}{2} + \tau(1 - 2mw)\enspace.
    \end{align*}
    For the right hand side $\Exp_\mu (\eta_\sigma(X) - 1/4)_+$, there are two cases $\tau > 1/4$ and $0< \tau \leq 1/4$, one can easily show that as long as $b' \leq 1/8$ there are no values of $\tau$ which allow to fix $\theta^* = 1/4$.
    Therefore, $\tau > 1/4$ and we can write for every $\sigma \in \{-1, 1\}$
    \begin{align*}
        \Exp_\mu (\eta_\sigma(X) - 1/4)_+ &= \sum_{j = 1}^m\int_{\class{X}_j} (                                      \sigma_j\xi(x))_+ d\mu(x)
                                                +
                                                \sum_{j = 1}^m\int_{\class{X}_{-j}} (-\sigma_j\xi(x))_+ d\mu(x)
                                        +\int_{A_0}(\tau - 1/4) d\mu(x)\\
                                        &= mwb' + (\tau - 1/4)(1 - 2mw)\enspace.
    \end{align*}
    Finally, the parameter $\tau$ must satisfy the following equality
    \begin{align*}
        \frac{1}{4}\left(\frac{mw}{2} + \tau(1 - 2mw) \right) = mwb' + (\tau - 1/4)(1 - 2mw)\enspace,
    \end{align*}
    solving for $\tau$ we get
    \begin{align*}
        \tau = \frac{1}{3} + \left(\frac{1}{12} - \frac{2b'}{3}\right)\left(\frac{2mw}{1 - 2mw}\right)\enspace.
    \end{align*}
    If $mw \leq 1/2$ we can ensure that the value of $\tau \leq 1$, that is, it is a valid choice for the regression function.
    Let us demonstrate that the margin assumption~\ref{ass:margin_assumption} holds for an appropriate choice of $m$ and $w$.
    Define $x_0 = (1/2q, \ldots, 1/2q)^\top$, then for every $\sigma \in \{-1, 1\}$ we have
    \begin{align*}
        P_X(0 < \abs{\eta_\sigma(X) - 1/4} \leq \delta)
        &=
        \frac{2mw}{\lambda(\mathcal{B}(0, (4q)^{-1}))}\int_{\mathcal{B}(x_0, (4q)^{-1})}\ind{C_\varphi q^{-\beta} u(q\norm{x - n_q(x)}_2) \leq \delta} dx\\
        &\phantom{=}
        +
        \frac{1 - 2mw}{\lambda(A_0)}\int_{A_0}\ind{\frac{1}{3} + \left(\frac{1}{12} - \frac{2b'}{3}\right)\left(\frac{2mw}{1 - 2mw}\right) - \frac{1}{4} \leq \delta} dx\\
        &= 2mw\ind{\delta \geq C_{\varphi}q^{-\beta}} + \frac{1 - 2mw}{\lambda(A_0)}\int_{A_0}\ind{\frac{1}{12} + \left(\frac{1}{12} - \frac{2b'}{3}\right)\left(\frac{2mw}{1 - 2mw}\right) \leq \delta} dx\enspace,
    \end{align*}
    as long as $b' \leq 3 /24$ we can continue as
    \begin{align*}
        P_X(0 < \abs{\eta_\sigma(X) - 1/4} \leq \delta)
        &\leq
        2mw\ind{\delta \geq C_{\varphi}q^{-\beta}} + \ind{\delta \geq \frac{1}{12}}\\
        &\leq 2mw\ind{\delta \geq C_{\varphi}q^{-\beta}} + 12^{\alpha}\delta^{\alpha}\enspace.
    \end{align*}
    Therefore, if $mw$ is of order $q^{-\alpha \beta}$ the margin assumption is satisfied with $\delta_0 = 1/12$.
    % One can modify the construction above to match any fixed constant $\delta_0 \in (1/12, 1/2)$ (for $\delta \in (0, 1/12)$ the above construction yields the desired result), for this it is sufficient to set $\tau = \theta^* + 1/2$ and find the value of the optimal threshold $\theta^*$ that satisfies the fixed point equation.
    % Though, the modification is straightforward, it involves a lot of standard algebraic manipulations and does not provide with extra insights.
    % That is way we omit these details to facilitate the understanding.
    The strong density assumption can be checked similarly to~\cite{Audibert_Tsybakov07}.
    To finish the prove, for every $\sigma \in \{-1, 1\}^m$ we denote by $P^\sigma$ the distribution of $(X, Y)$ with the marginal $P_X$ and the regression function $\eta^\sigma$.
    Thus, one can write for any $\hat g$
    \begin{align*}
        \sup_{\Prob \in \class{P}_{\Sigma}} \Exp_{(\data_n, \data_N)}\excess(\hat g)
        \geq
        \sup_{\sigma \in \{-1, 1\}^m} \frac{1}{2}\Exp^{\sigma}_{(\data_n, \data_N)}\sum_{i = -m, i\neq 0}^m\Exp_{P_X}\abs{\varphi(X)}\ind{(1 + \sign(i)\sigma_i)/2 \neq \hat g(X)}\ind{X \in \class{X}_i}\enspace,
    \end{align*}
    where $\Exp^{\sigma}_{(\data_n, \data_N)}$ is the expectation taken \wrt to the \iid realizations of $\data_n$ and $\data_N$ from $P^{\sigma}$ and $P_X$ respectively, and $\sign(i) = 1$ if $i > 0$ and $\sign(i) = -1$ if $i < 0$.
    The rest of the proof is obtained following the proof of~\cite[Lemma 5.1.]{Audibert04} and in particular the chain of inequalities in~\cite[Eq.~(6.26)]{Audibert04}.
    That is, we get for some $C > 0$ independent from $N, n$
    \begin{align*}
        \sup_{\Prob \in \class{P}_{\Sigma}} \Exp_{(\data_n, \data_N)}\excess(\hat g) \geq C mwq^{-\beta}(1 - C_{\varphi}q^{-\beta}\sqrt{n w})
    \end{align*}
    Finally, we conclude by setting the parameters $m, w, q$ as
    \begin{align*}
        q = \floor{\bar C n^{\frac{1}{2\beta + d}}},\quad w = C' q^{-d},\quad m = \floor{C'' q^{d - \alpha \beta}}\enspace.
    \end{align*}
    Note that thanks to the condition $\alpha \beta \leq d$ such a choice is is always valid for appropriately chosen constants $\bar C, C', C''$.
\end{proof}

\bibliographystyle{spmpsci}      % mathematics and physical sciences
\bibliography{references_all}   % name your BibTeX data base

\providecommand{\AC}{A.-C}\providecommand{\CA}{C.-A}\providecommand{\CH}{C.-H}\providecommand{\CJ}{C.-J}\providecommand{\JC}{J.-C}\providecommand{\JP}{J.-P}\providecommand{\JB}{J.-B}\providecommand{\JF}{J.-F}\providecommand{\JM}{J.-M}\providecommand{\KW}{K.-W}\providecommand{\PL}{P.-L}\providecommand{\RE}{R.-E}\providecommand{\SJ}{S.-J}\providecommand{\XR}{X.-R}\providecommand{\WX}{W.-X}
\begin{thebibliography}{10}
\providecommand{\url}[1]{{#1}}
\providecommand{\urlprefix}{URL }
\expandafter\ifx\csname urlstyle\endcsname\relax
  \providecommand{\doi}[1]{DOI~\discretionary{}{}{}#1}\else
  \providecommand{\doi}{DOI~\discretionary{}{}{}\begingroup
  \urlstyle{rm}\Url}\fi

\bibitem{Audibert04}
Audibert, J.Y.: Aggregated estimators and empirical complexity for least square
  regression.
\newblock Ann. Inst. H. Poincar{\'e} Probab. Statist. \textbf{40}(6), 685--736
  (2004)

\bibitem{Audibert_Tsybakov07}
Audibert, J.Y., Tsybakov, A.B.: Fast learning rates for plug-in classifiers.
\newblock Ann. Statist. \textbf{35}(2), 608--633 (2007)

\bibitem{Bartlett_Mendelson02}
Bartlett, P.L., Mendelson, S.: Rademacher and {G}aussian complexities: risk
  bounds and structural results.
\newblock J. Mach. Learn. Res. \textbf{3}(Spec. Issue Comput. Learn. Theory),
  463--482 (2002)

\bibitem{Bobkov_Ledoux16}
Bobkov, S., Ledoux, M.: One-dimensional empirical measures, order statistics
  and {Kantorovich} transport distances  (2016).
\newblock To appear in the {Memoirs of the Amer. Math. Soc.}

\bibitem{Chzhen_Denis_Hebiri19}
Chzhen, E., Denis, C., Hebiri, M.: {Minimax semi-supervised confidence sets for
  multi-class classification} (2019).
\newblock Preprint, \url{https://arxiv.org/abs/1904.12527}

\bibitem{Conte_Boor80}
Conte, S., Boor, C.: Elementary Numerical Analysis: An Algorithmic Approach,
  3rd edn.
\newblock McGraw-Hill Higher Education (1980)

\bibitem{Dembczynski_Kotlowski_Koyejo_Natarajan17}
Dembczynski, K., Kot{\l}owski, W., Koyejo, O., Natarajan, N.: Consistency
  analysis for binary classification revisited.
\newblock In: ICML, pp. 961--969. JMLR. org (2017)

\bibitem{Dvoretzky_Kiefer_Wolfowitz56}
Dvoretzky, A., Kiefer, J., Wolfowitz, J.: Asymptotic minimax character of the
  sample distribution function and of the classical multinomial estimator.
\newblock Ann. Math. Statist. \textbf{27}(3), 642--669 (1956)

\bibitem{Keerthi_Sindhwani_Chapelle07}
Keerthi, S., Sindhwani, V., Chapelle, O.: An efficient method for
  gradient-based adaptation of hyperparameters in svm models.
\newblock In: NIPS, pp. 673--680 (2007)

\bibitem{Koyejo_Natarajan_Ravikumar_Dhillon14}
Koyejo, O., Natarajan, N., Ravikumar, P., Dhillon, I.: Consistent binary
  classification with generalized performance metrics.
\newblock In: NIPS, pp. 2744--2752 (2014)

\bibitem{Lewis95}
Lewis, D.: Evaluating and optimizing autonomous text classification systems.
\newblock In: ACM, pp. 246--254. ACM Press (1995)

\bibitem{Massart90}
Massart, P.: The tight constant in the dvoretzky-kiefer-wolfowitz inequality.
\newblock Ann. Probab. \textbf{18}(3), 1269--1283 (1990)

\bibitem{Massart_Nedelec06}
Massart, P., N{\'e}d{\'e}lec, {\'E}.: Risk bounds for statistical learning.
\newblock Ann. Statist. \textbf{34}(5), 2326--2366 (2006)

\bibitem{Menon_Narasimhan_Agarwal_Chawla13}
Menon, A., Narasimhan, H., Agarwal, S., Chawla, S.: On the statistical
  consistency of algorithms for binary classification under class imbalance.
\newblock In: ICML, vol.~28, pp. 603--611. PMLR (2013)

\bibitem{Narasimhan_Vaish_Agarwal14}
Narasimhan, H., Vaish, R., Agarwal, S.: On the statistical consistency of
  plug-in classifiers for non-decomposable performance measures.
\newblock In: NIPS, pp. 1493--1501 (2014)

\bibitem{Rigollet07}
Rigollet, P.: Generalization error bounds in semi-supervised classification
  under the cluster assumption.
\newblock Journal of Machine Learning Research \textbf{8}(Jul), 1369--1392
  (2007)

\bibitem{Rigollet_Vert09}
Rigollet, P., Vert, R.: Optimal rates for plug-in estimators of density level
  sets.
\newblock Bernoulli  (2009)

\bibitem{Singh_Nowak_Zhu09}
Singh, A., Nowak, R., Zhu, J.: Unlabeled data: Now it helps, now it doesn{'}t.
\newblock In: NIPS, pp. 1513--1520 (2009)

\bibitem{Tsybakov09}
Tsybakov, A.B.: Introduction to nonparametric estimation.
\newblock Springer Series in Statistics. Springer, New York (2009)

\bibitem{Vallender74}
Vallender, S.: Calculation of the wasserstein distance between probability
  distributions on the line.
\newblock Theory of Probability \& Its Applications \textbf{18}(4), 784--786
  (1974)

\bibitem{Rijsbergen74}
{van Rijsbergen}, C.: Foundation of evaluation.
\newblock Journal of documentation \textbf{30}(4), 365--373 (1974)

\bibitem{Vapnik98}
Vapnik, V.N.: Statistical learning theory.
\newblock Wiley (1998)

\bibitem{Yan_Koyejo_Zhong_Ravikumar18}
Yan, B., Koyejo, S., Zhong, K., Ravikumar, P.: Binary classification with
  karmic, threshold-quasi-concave metrics.
\newblock In: ICML, vol.~80. PMLR (2018)

\bibitem{Yang99}
Yang, Y.: Minimax nonparametric classification: Rates of convergence.
\newblock IEEE Transactions on Information Theory \textbf{45}(7), 2271--2284
  (1999)

\bibitem{Ye_Chai_Lee_Chieu12}
Ye, N., Chai, K., Lee, W., Chieu, H.: Optimizing f-measures: A tale of two
  approaches.
\newblock In: ICML (2012)

\bibitem{Zhao_Edakunni_Pocock_Brown13}
Zhao, M.J., Edakunni, N., Pocock, A., Brown, G.: Beyond fano's inequality:
  bounds on the optimal f-score, ber, and cost-sensitive risk and their
  implications.
\newblock JMLR \textbf{14}(Apr), 1033--1090 (2013)

\end{thebibliography}

% \begin{thebibliography}{}
% %
% % and use \bibitem to create references. Consult the Instructions
% % for authors for reference list style.
% % Format for Journal Reference
% Author, Article title, Journal, Volume, page numbers (year)
% Author, Book title, page numbers. Publisher, place (year)
% \end{thebibliography}

\end{document}